\documentclass[12pt,a4paper,reqno]{amsart}

\newtheorem{theorem}{Theorem}
\newtheorem{lemma}{Lemma}

\newcommand{\gauss}[2]{\genfrac{[}{]}{0pt}{}{#1}{#2}_q}  
\newcommand{\gausss}[2]{\genfrac{[}{]}{0pt}{}{#1}{#2}_{1/q}}

\usepackage{amsmath}
\usepackage{amssymb}
\usepackage{amsfonts}

\advance\textheight by \topskip
\title[A $q$--analogue of a formula of Hernandez]
{A $q$--analogue of a formula of Hernandez obtained
by inverting a result of Dilcher}

\author{Helmut Prodinger}
\thanks{This research was conducted while the author was
a host of the projet Algo at INRIA, Rocquencourt. The funding
came from the Austrian--French ``Amad\'ee'' cooperation.}

\address{ Helmut Prodinger,
Centre for Applicable Analysis and Number Theory,
 Department of Mathematics,
University of the Witwatersrand, P.~O. Wits, 
2050 Johannesburg, South Africa, email:
{\tt helmut@gauss.cam.wits.ac.za}.
}

\date{July 7, 1999}

\begin{document}
\begin{abstract}
We prove a $q$--analogue of the formula
\begin{equation*}
\sum_{1\le k\le n} \binom nk(-1)^{k-1}\sum_{1\le i_1\le i_2\le\dots \le
i_m=k}\frac1{i_1i_2\dots i_m}= \sum_{1\le k\le n}\frac{1}{k^m} 
\end{equation*}
by inverting a formula due to Dilcher.
\end{abstract}

\maketitle

\section{The identities}

Hern\'andez in \cite{Hernandez99} proved the following identity:

\begin{equation}\label{Hernandez}
\sum_{1\le k\le n} \binom nk(-1)^{k-1}\sum_{1\le i_1\le i_2\le\dots \le
i_m=k}\frac1{i_1i_2\dots i_m}= \sum_{1\le k\le n}\frac{1}{k^m} .
\end{equation}
However this identity does not really require a proof, since we will
show that it is just an inverted form of an identity of Dilcher;
\cite{Dilcher95}
\begin{equation}\label{dilcher}
\sum_{1\le k\le n} \binom nk(-1)^{k-1}\frac1{k^m}
=\sum_{1\le i_1\le i_2\le\dots \le
i_m\le n}\frac1{i_1i_2\dots i_m}.
\end{equation}

Define for $k\ge1$
$$
a_k:=-\sum_{1\le i_1\le i_2\le\dots \le
i_m=k}\frac1{i_1i_2\dots i_m}
\qquad\text{and}\qquad
b_k:=\frac1{k^m},
$$
and $a_0=b_0=0$, then the identities are

\begin{align*}
\sum_{0\le k\le n}\binom nk(-1)^ka_k&=\sum_{0\le k\le n}b_k,\\
\sum_{0\le k\le n}\binom nk(-1)^kb_k&=\sum_{0\le k\le n}a_k.
\end{align*}
They are {\sl inverse relations,} as can be seen by introducing ordinary
generating functions $A(z)=\sum a_nz^n$ and
$B(z)=\sum b_nz^n$. Then they are

\begin{align*}
\frac1{1-z}A\Big(\frac z{z-1}\Big)&=\frac1{1-z}B(z),\\
\frac1{1-z}B\Big(\frac z{z-1}\Big)&=\frac1{1-z}A( z).
\end{align*}
However

$$
w=\frac z{z-1} \longleftrightarrow z=\frac w{w-1},
$$
and the proof is finished

We note that
Dilcher's sum appears also in disguised form in \cite{FlSe95}.

\section{A $q$--analogue}

Dilcher's formula (\ref{dilcher}) was only a corollary of his elegant
 $q$--version;

$$
\sum_{1\le k\le n} \gauss
nk(-1)^{k-1}\frac{q^{\binom{k+1}{2}+(m-1)k}}{(1-q^k)^m} =\sum_{1\le i_1\le
i_2\le\dots \le i_m\le n}\frac{q^{i_1}}{1-q^{i_1}}\dots \frac{
q^{i_m}}{1-q^{i_m}}. $$
Here, $\gauss nk$ denotes the Gaussian polynomial
\begin{equation*}
\gauss nk
=\frac{(q;q)_n}{(q;q)_k(q;q)_{n-k}}
\end{equation*}

with
\begin{equation*}
(x;q)_n:=(1-x)(1-xq)\dots(1-xq^{n-1}).
\end{equation*}
Apart from Dilcher's paper \cite{Dilcher95}, the article
\cite{AnCrSi97} is also of some relevance in this context.

Therefore it is a natural 
question to find a 
$q$--analogue of Hern\'andez' formula, or, what amount to the same,
to find the appropriate inverse relations for the $q$--analogues.

We state them in the following lemma which is almost surely
not new. However, I found it more appealing to derive it myself
rather than hunting for it in the vast $q$--literature.

\begin{lemma}
\begin{align}
\sum_{k=0}^nb_k&=\sum_{k=0}^n\gauss nk(-1)^kq^{\binom k2}a_k, \label{first}\\
\sum_{k=0}^nq^{-k}a_k&=\sum_{k=0}^n\gauss nk(-1)^kq^{-kn+\binom k2}b_k \label{second}.
\end{align}
\end{lemma}

{\bf Proof}. First note that it is much easier to prove this than to find
it (it took me a few hours, but I am not too experienced).

We note that, as is usual in such cases, it is sufficient to
prove it for a basis of the vector space of polynomials. Here, we choose
$a_n=q^nx^n(1-\frac1x)$ for $n\ge1$ and $a_0=1$.

We need the following standard formul{\ae} that are consequences
of the $q$--binomial theorem (see e.~g. \cite{Andrews76}):
\begin{align*}
\sum_{k=0}^n
\gauss nk(-1)^kx^k&=
(x;q)_n,\\
\sum_{k=0}^n
\gauss nk(x;q)_{k}x^{n-k}&=1
\end{align*}

We plug the form of $a_n$ into the right hand side of (\ref{first})
and obtain
\begin{align*}
(1-\tfrac1x)\sum_{k=0}^n
\gauss nk(-1)^kq^kx^k=
(1-\tfrac1x)(qx;q)_n=-\tfrac1x(x;q)_{n+1}=
\sum_{k=0}^nb_k.
\end{align*}

Thus
\begin{align*}
b_n=-\tfrac1x(x;q)_{n+1}+\tfrac1x(x;q)_{n}=
-\tfrac1x(x;q)_{n}(1-xq^n-1)=q^n(x;q)_n.
\end{align*}
We are done if these values of $a_n$ and $b_n$ also satisfy
the relation (\ref{second}). Note that (\ref{second})
can be rewritten as

\begin{equation*}
\sum_{k=0}^nq^{-k}a_k=\sum_{k=0}^n\gausss nk(-1)^k
q^{-\binom {k+1}2}b_k .
\end{equation*}

We plug  $b_n$ into the right hand side of (\ref{second})
and obtain

\begin{align*}
\sum_{k=0}^n\gausss nk(-1)^k
q^{-\binom {k+1}2}q^k(x;q)_k=
\sum_{k=0}^n\gausss nk(-1)^k
x^k(\tfrac 1x;\tfrac 1q)_k=x^n=
\sum_{k=0}^nq^{-k}a_k.
\end{align*}
Thus
\begin{align*}
q^{-n}a_n=x^n-x^{n-1}=x^n(1-\tfrac1x).
\end{align*}

We would like to remark that an alternative proof could be found by
dealing with the matrices of {\sl connecting coefficients.}

Define matrices

\begin{align*}
T:=\left[\gauss nk(-1)^kq^{\binom {k}2}\right]_{n,k},
\qquad
S:=\left[\gauss nk(-1)^kq^{-kn+\binom {k}2}\right]_{n,k},
\end{align*}

$$U=[\mathbf 1_{n\ge k}]_{n,k},
\quad\text{and}\quad
V=[q^{-k}\mathbf 1_{n\ge k}]_{n,k}.
$$
Then we have to prove that $S=VT^{-1}U$.

This is not too hard, since
\begin{align*}
T^{-1}=\left[\gauss nk(-1)^kq^{-kn+\binom {k+1}2}\right]_{n,k},
\end{align*}
and
\begin{align*}
T^{-1}U=\left[\gauss {n-1}{k-1}(-1)^kq^{-n(k-1)+\binom {k}2}\right]_{n,k}.
\end{align*}

\begin{theorem}{\rm [$q$--analogue of Hern\'andez' formula]}
\begin{align*}
\sum_{1\le k \le n}\gauss nk (-1)^{k-1}q^{-kn+\binom k2}
\sum_{1\le i_1\le
i_2\le\dots \le i_m=k}\frac{q^{i_1}}{1-q^{i_1}}\dots \frac{
q^{i_m}}{1-q^{i_m}}=
\sum_{1\le k \le n}\frac{q^{k(m-1)}}{(1-q^k)^m}.
\end{align*}

\end{theorem}

\bibliographystyle{plain}


\end{document}